\documentclass[11pt, a4paper]{amsart} 
\usepackage{graphicx,psfrag} 

\newcommand{\inv}{^{-1}}

\newcommand{\tto}{\buildrel * \over \rightarrow }

\newcommand{\abs}[1]{\vert #1\vert} 

\newcommand{\set}[1]{\{ #1 \}} 
\newcommand{\subgroup}[1]{\langle #1\rangle} 
\newcommand{\ovr}[1]{\overline{ #1 }} 
\newcommand{\be}{\begin{enumerate}}
\newcommand{\ee}{\end{enumerate}}

\newcommand{\G}{\Gamma} 
\newcommand{\Si}{\Sigma}

\newcommand{\e}{\varepsilon}

\newcommand{\mt}{\emptyset}

\newtheorem{theorem}{Theorem}[section]
\newtheorem{lemma}[theorem]{Lemma}
\newtheorem{corollary}[theorem]{Corollary} 
\theoremstyle{definition}
\newtheorem{definition}[theorem]{Definition}
\newtheorem{example}[theorem]{Example} 
\newtheorem{question}{Question}
\theoremstyle{remark}

\title{Word Hyperbolic Semigroups}

\author{Andrew Duncan} \address{School of Mathematics and Statistics,
University of Newcastle upon Tyne} \email{A.Duncan@ncl.ac.uk}

\author{Robert H. Gilman} \address{Department of Mathematical
Sciences, Stevens Institute of Technology, Hoboken, New Jersey 07030}
\email{rgilman@stevens-tech.edu}

\begin{document}

\maketitle

\section{Introduction}

The study of automatic and word hyperbolic groups is a prominent topic in
geometric group theory.  The definition of automatic groups extends
naturally to semigroups~\cite{DR, OSM}, but word hyperbolic groups are
defined by a geometric condition, the thin triangle condition, whose
extension is not immediate.  The problem is that for semigroups the
metric properties of the Cayley diagram are not as closely tied to the
algebraic structure as they are in the case of groups.  We propose a
language--theoretic definition of word hyperbolic semigroup which for
groups is equivalent to the original definition.

We also consider briefly word problems of semigroups. The standard
word problem of a group is the language of all words representing the
identity.  As it stands this definition makes no sense for semigroups,
but there is a variation which does.

\section{Background}\label{background}

We assume familiarity with formal language theory;
\cite[Volume~1]{handbook} affords a good introduction.  A brief
discussion is included here for the reader's convenience and to
establish notation.

\subsection{Rational Subsets}\label{rationalsets}

The rational subsets of a monoid are the closure of its finite
subsets under union, product, and generation of submonoid. It follows
from this definition that non--empty rational sets are determined by
expressions built up from the elements of the monoid together with $+$ for
union, $^*$ for generation of submonoid, and concatenation for product. For
example $m_1+m_2m_3^*$ denotes $\set{m_1, m_2m_3^i \mid i\ge 0}$. We
ignore the distinction between expressions and sets and write
$m_1+m_2m_3^*=\set{m_1, m_2m_3^i \mid i\ge 0}$ etc. Also we write
$\mt$ for the empty set and $L^+$ for $LL^*$.

Rational subsets of $\Sigma^*$, the free monoid of all words on a
finite alphabet $\Sigma$, have a special name. They are called regular
languages over $\Sigma$. Regular languages are closed under
intersection and complement. Rational subsets in general are not.  

\subsection{Automata}\label{automata}

Rational subsets of monoids may be defined as the sets accepted by
finite automata.  A 
finite automaton $\mathcal A$ over a monoid $M$ is a finite directed
graph with edge labels from $M$, a distinguished initial vertex, and
some distinguished terminal vertices. The subset of $M$ accepted by
$\mathcal A$ is the set of labels of directed paths beginning at the initial
vertex and ending at one of the terminal vertices. The label of a path
is the product of its edge labels in order. Paths of length zero have
label $1$, the identity element of $M$. 

\begin{figure}[hbt] \centering 
\psfrag{e}[cc]{$\e$}
\psfrag{i}[cc]{$a$}
\psfrag{j}[cc]{$b$} 
\includegraphics{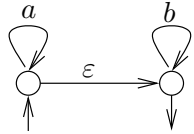}
\caption{An automaton accepting $a^*b^*$.\label{fa}}
\end{figure}

From now on we say automaton instead of finite
automaton. Figure~\ref{fa} shows an 
automaton accepting the 
%rational language 
regular language %ajd  
$a^*b^*= \set{a^ib^j \vert i,j,\ge
0}$. The symbol $\e$ in that figure denotes the empty word, the
identity element of $\Sigma^*$. The initial vertex is shown as the
target of an edge without any source, and similarly for the single
terminal vertex. 

Finite automata over $\Sigma^*\times M$ accept rational subsets of
$\Sigma^*\times M$, but they may also be used to define languages over
$\Sigma$. (We shall use this technique only in
Example~\ref{bicyclic}.) The language defined by an automaton
$\mathcal A$ over $\Sigma^*\times M$ is the set of words $w$ such that
$(w,1)$ is accepted by $\mathcal A$.

\begin{figure}[hbt] 
\centering 
\psfrag{e}[cc]{$(\epsilon,1)$}
\psfrag{i}[cc]{$(a,p_1)$}
\psfrag{j}[cc]{$(b,q_1)$} 
\includegraphics{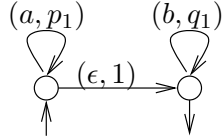}
\caption{An automaton defining $\{a^ib^i \mid i\ge 0\}$. \label{pda0}}
\end{figure}

For example let $M_{cf}$ be the monoid presented by $\langle p_i, q_i
\mid p_iq_i=1, i>0  
\rangle$. Figure~\ref{pda0} shows an automaton over $\Sigma^*\times
M_{cf}$ which accepts the rational set 
$R=(a,p_1)^*(b,q_1)^*=\{(a^ib^j,p_1^iq_1^j)\}$ and defines the
language $\set{a^ib^i}$. 
In~\cite{RG1} it is shown that automata over $\Sigma^*\times M_{cf}$ define the
class of context--free languages over $\Sigma$.

\subsection{Rational Transductions}

A rational transduction $\rho:\Sigma^*\to \Delta^*$ from one finitely
generated free monoid to another is a rational binary relation, that
is, a rational subset of $\Sigma^*\times\Delta^*$. 
It follows that rational transductions are closed
under union, product and generation of submonoids. They are also
closed under inverse and composition in the sense of binary
relations. Images of regular and context--free languages under
rational transductions are regular and context--free respectively. In
particular regular and context--free languages are closed under
homomorphism, inverse homomorphism, and intersection with regular
languages.

We need one more result. For any word $w$ in a free monoid,
$w^r$ is $w$ written from right to left. Note that $(wv)^r=v^rw^r$.

\begin{lemma}\label{reverse-lemma} For any homomorphism
$h:\Sigma^+\to\Delta^+$ between finitely generated free semigroups
the map $\tau_h:\Sigma^+\to\Delta^+$ defined by $\tau_h(x)= h(x^r)^r$ is a 
rational transduction.
\end{lemma}
\begin{proof}
The graph of $\tau_h$ is $\bigl(\sum_{a\in \Sigma}(a,h(a)^r)\bigr)^*$. 
\end{proof}

\subsection{Context--Free Grammars}

A context--free grammar $\mathcal G$ consists of a collection of
disjoint finite
sets. Namely a terminal alphabet, $\Sigma$, a set of nonterminals, $N$,
and a set of productions of the form $A\to\alpha$ where $A\in N$ and
$\alpha$ is 
%a sentential form, that is %ajd 
a word in $(\Sigma+N)^*$.
There is also a designated start symbol $S\in N$.

For 
%sentential forms %ajd
$\alpha$ and $\beta$ 
in $(\Sigma+N)^*$ %ajd
the notation
$\alpha\to \beta$ means that the lefthand side of some production is a
subword of $\alpha$, and that $\beta$ is obtained by replacing that
subword by the righthand side of that production. The effect of zero or
more replacements is denoted by $\alpha \tto \beta$. When
$\alpha\tto\beta$, we say $\beta$ is derived from $\alpha$ or $\alpha$
derives $\beta$. The language generated by the grammar $\mathcal G$ is
$\set{w\mid w\in \Sigma^*, S\tto w}$.

Context--free grammars generate context--free
languages. Every context--free language not containing the empty word
%is 
can be  %ajd
generated by a context--free grammar in Chomsky normal form,
that is a grammar with all productions of the form $A\to BC$ or $A\to
a$ for nonterminals  
%$A, B$ and terminals $a$. 
$B, C$ and terminal $a$. Context--free languages are closed under 
finite union and under 
intersection with regular languages. %ajd
\section{Definitions and Initial Results}

\begin{definition} A choice of generators for a semigroup $S$ is a
  surjective homomorphism $f:\Sigma^+\to S$ from the free semigroup
  $\Sigma^+$ over a finite alphabet $\Sigma$. We write $\ovr w$ for
  $f(w)$ when $w\in\Sigma^*$, and $\ovr X$ for $f(X)$ if
  $X\subset\Sigma^*$.    
\end{definition}

\begin{definition} A combing with respect to a choice of generators
  $\Sigma^+\to S$ is a subset $R\subset \Sigma^+$ such that $\ovr R =
  S$. A regular combing is a combing which is a regular language.
\end{definition}

In order to state our definition of hyperbolic monoid we need one more
bit of notation. For any alphabet $\Sigma$, $\#$ denotes a symbol not
in $\Sigma$; and (in the notation of Section~\ref{background})
$\Sigma_\#=\Sigma+\#$. 

\begin{definition}\label{hyperbolic} A semigroup is word
  hyperbolic if for some choice of generators $f:\Sigma^+\to S$ and regular
  combing $R\subset \Sigma^+$ the language $T=\{u\#v\#w^r \mid u,v,w\in
  R\text{ and }\ovr u \ovr v = \ovr w\}$ is context--free. 
\end{definition}

We call $T$ the multiplication table determined by $R$ and refer
to $R$ and $T$ collectively as a hyperbolic structure for $M$. From now
on we abbreviate word hyperbolic by hyperbolic and regular
combing by combing. 

For groups there are two notions of
hyperbolicity. A group $G$ may be hyperbolic according to
Definition~\ref{hyperbolic} or as in the original geometric
definition. We say that $G$ is hyperbolic in the first case and
hyperbolic as a group in the second.

\begin{theorem}\label{generators} If $S$ is hyperbolic with respect to one
  choice of generators, then it is hyperbolic with respect to every
  choice.
\end{theorem} 

\begin{proof} Let $S$ be hyperbolic for the choice of generators
  $f:\Sigma^+\to S$, and suppose $g:\Delta^+\to S$ is a second 
  choice of generators.  Pick a homomorphism $h:\Sigma^+\to \Delta^+$ such that
  $g\circ h = f$.  $R_1=h(R)$ is a regular combing; let $T_1$ be the
  corresponding multiplication table. To complete the
  proof it suffices to exhibit a rational transduction
  $\rho:\Sigma_\#^*\to \Delta_\#^*$ with $\rho(T)=T_1$, for then we
  know that $T_1$ is context--free. 

  By Lemma~\ref{reverse-lemma} the product of
  transductions $\rho=(h)(\#,\#)(h)(\#,\#)(\tau_h)$ is rational.
  It is easy to see that $\rho$ is a partial function with domain
  $\Sigma^+\#\Sigma^+\#\Sigma^+$ and 
  $\rho(u\#v\#w^r)=h(u)\#h(v)\#h(w)^r$. It follows by a straightforward
  argument that $\rho(T)=T_1$.
\end{proof}

For a monoid $M$ it is natural to consider a choice of generators to
be a surjective homomorphism $\Sigma^*\to M$ from the free
monoid $\Sigma^*$ over $\Sigma$.
   
\begin{theorem}\label{monoid} Let $M$ be a monoid and $f:\Sigma^*\to
  M$ a surjective homomorphism. $M$ is hyperbolic if and only if there
  is a regular 
  combing $R\subset \Sigma^*$ such that $T=\{u\#v\#w^r \mid u,v,w\in
  R\text{ and }\ovr u \ovr v = \ovr w\}$ is context--free.
\end{theorem}

\begin{proof}
  Suppose $f:\Sigma^*\to  M$ is any surjective homomorphism.
  Let $\Delta= \Sigma+x$ where $x$ is a letter not in $\Sigma$. Define
  $h:\Delta^+\to \Sigma^*$ by $h(x)=\e$ and $h(a)=a$ for $a\in \Sigma$.
  Extend $h$ to $h:\Delta_\#^+\to \Sigma_\#^*$ by $h(\#)=\#$. Note
  that as images of letters in $\Delta_\#$ have length at most $1$ in
  $\Sigma_\#^*$, $h(z^r)=h(z)^r$ for any $z\in\Delta_\#^+$. 

  The composition $f\circ h:\Delta^+\to M$ is a semigroup choice of
  generators for $M$. If $M$ is hyperbolic, then by
  Theorem~\ref{generators} there is a regular combing $R_1\subset
  \Delta^+$ whose associated multiplication table $T_1$ is
  context--free. Let $R=h(R_1)$, and check that $T=h(T_1)$ is the
  multiplication table for $R$. Conversely if $R$ and $T$ are as in
  the statement of the theorem, observe that $R_1=h\inv(R)$ is a
  regular combing with corresponding multiplication table
  $T_1=h\inv(T)$.
\end{proof}

For any semigroup $S$ the semigroup obtained by
adjoining an identity element, $1$, to $S$ is denoted $S^I$. Likewise the
semigroup obtained by adjoining a zero, $0$, is $S^Z$.

\begin{theorem}\label{extend} $S^Z$ is hyperbolic if and only if $S$
is hyperbolic. If $S^I$ is hyperbolic, then so is $S$.
\end{theorem}

\begin{proof}
  Suppose that $S$ is hyperbolic with respect to the choice of generators
  $\Sigma^+\to S$ and combing $R\subset \Sigma^+$. Add a new symbol
  $x$ to $\Sigma$. Define $\ovr x = 0$, and include $x$ in
  $R$. The effect of these changes is to replace the
  multiplication table $T$ by $T+R\#x\#x+ x\#R\#x+x\#x\#x$, which is
  still context--free. Hence $S^Z$ is hyperbolic.

  To prove the converse assume $S^Z$ is hyperbolic. By
  Theorem~\ref{generators} generators $(\Sigma + x)^+\to S^Z$ chosen
  as in the preceding paragraph afford a hyperbolic structure $R_1$,
  $T_1$.  $R=R_1\cap \Sigma^+$ and
  $T=T_1\cap(\Sigma^+\#\Sigma^+\#\Sigma^+)$ are a hyperbolic structure
  for $S$.

  The proof that $S$ is hyperbolic if $S^I$ is proceeds
  similarly. Take a hyperbolic structure $R_1, T_1$ with respect to a
  choice of generators $(\Sigma + x)^+ \to S^I$ such that $\Sigma^+\to
  S$
  is a choice of generators for $S$ and $\ovr x = 1$. Construct a
  combing $R\subset \Sigma^+$ for $S$ by removing all occurrences of
  $x$ from words in $R_1$ and then discarding the empty word. We leave
  it to the reader to check that 
  %$R_1$  
  $R$ %ajd
  is regular and its
  multiplication table is context--free.
\end{proof}

\begin{question} Is $S^I$ hyperbolic if $S$ is? 
\end{question}

\noindent
A hyperbolic structure $R,T$ for $S$ extends to a hyperbolic
structure for $S^I$ if and only if 
$L=\set{u\#w^r\mid u,w\in R\text{ and }\ovr u= \ovr w}$ is context--free.  
However, the condition that $L$ be context--free does not
seem to follow from the assumption that $S$ is hyperbolic except in
the case that $R$ projects bijectively to $M$. For then
$L=\set{u\#u^r\mid u\in R}$, which is context--free for any
regular language $R$. Thus an affirmative answer to the next question
implies an affirmative answer to the previous one.

\begin{question} Does every hyperbolic monoid have a hyperbolic
  structure with uniqueness, that is, a hyperbolic structure $R$, $T$
  such that $R$ projects bijectively?
\end{question}

Now we give some examples of hyperbolic semigroups.

\begin{example}
  Finite semigroups are hyperbolic. If $S$ is a finite semigroup, take
  $R=\Si=S$. $R$ is a regular combing for $S$; $T$ is finite and hence
  context--free.
\end{example}

\begin{example}\label{bicyclic}
  The bicyclic monoid $M=\subgroup{a,b\mid ab=1}$ is hyperbolic.
  Take the choice of generators  $(a+b)^*\to M$
  determined by the preceding presentation. By Theorem~\ref{monoid} it
  suffices to show that the multiplication table $T$ corresponding to
  the combing $R=b^*a^*$ is context--free. We employ the method
  of Section~\ref{automata}.

  $T$ consists of all words $b^ia^j\#b^ka^l\#a^mb^n$ such that
  \begin{equation}
  j\ge k, j-k+l=m, i=n \text{ or } j< k, i+k-j=n, l=m. \label{eq1} 
  \end{equation}
  It suffices to
  check that $T$ is defined by the automaton $\mathcal A$ over
  $(a+b)^*\times M_{cf}$ given in Figure~\ref{pda}.

  \begin{figure}[hbt] 
  \centering 
  \psfrag{e}[cc]{$(\epsilon,1)$}
  \psfrag{h}[cc]{$(\#,1)$}
  \psfrag{i}[cc]{$(b,p_2)$}
  \psfrag{j}[cc]{$(a,p_1)$} 
  \psfrag{k1}[cc]{$(b,q_1)$}
  \psfrag{k2}[cc]{$(b,p_2)$} 
  \psfrag{l}[cc]{$(a,p_1)$}
  \psfrag{m}[cc]{$(a,q_1)$} 
  \psfrag{n}[cc]{$(b,q_2)$}
  \includegraphics{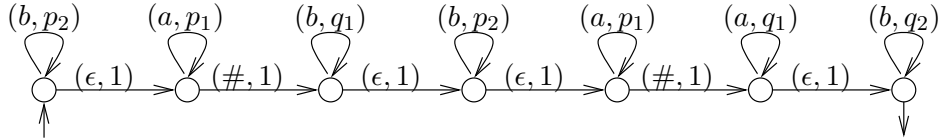}
  \caption{The automaton $\mathcal A$.\label{pda}}
  \end{figure} 

  Recall that $\mathcal A$
  defines the language of all words $w$ which occur in labels $(w,1)$
  of paths from the initial vertex of $\mathcal A$ to the final
  vertex. The labels of successful paths have the form
  $$
  (w,x)=(b^ia^j\#b^{k'+k''}a^l\#a^mb^n,
  p_2^ip_1^jq_1^{k'}p_2^{k''}p_1^lq_1^mq_2^n).
  $$ 
  Suppose $x=1$. Since $q_1$ has no right inverse, $j\ge k'$; and as
  $p_1$ has no left inverse, $l\le m$. Thus
  $x=p_2^ip_1^{j-k'}p_2^{k''}q_1^{m-l}q_2^n$. Either $k''=0, j-k'= m-l,
  i=n$ or 
 %$k''> 0, j=k', m=l, i+k'=n$. 
  $k''> 0, j=k', m=l, i+k^{''}=n$. %ajd
  In both cases $w$ satisfies
  Equation~\ref{eq1}.  
  
  Conversely every $w$ satisfying Equation~\ref{eq1} is in the
  language defined by $\mathcal A$. Set $k'=k, k''=0$ when $j\ge k$,
  and $k'=j$, $k''=k-j$ when $j<k$.
\end{example}

\begin{example}
Finitely generated free semigroups are hyperbolic. Take the
identity map $\Delta^+\to \Delta^+$ as a choice of generators.
$R=\Delta^+$ and 
%$T=\set{u\#v\#v^ru^r \mid u, v \in \Delta^*}$ 
 $T=\set{u\#v\#v^ru^r \mid u, v \in \Delta^+}$ %ajd
form a
hyperbolic structure.
\end{example}

\begin{example}\label{free}
  Finitely generated subsemigroups of free semigroups are
  hyperbolic. Observe that it suffices to check the assertion for
  subsemigroups of finitely generated free semigroups.  Suppose
  $h:\Si^+\to S\subset \Delta^+$ is a choice of generators for $S$,
  and let $T$ be the multiplication table of $\Delta^+$ from the
  preceding example. Check that the multiplication table determined by
  the combing $R=\Si^+$ of $S$ is $\rho^{-1}(T)$ where $\rho$ is the
  rational transduction from the proof of Theorem~\ref{generators}.
\end{example}

\begin{example}\label{subfree}
  Let $F$ be the free semigroup generated by $a$, $b$ and $c$, and let
  $S$ be the subsemigroup generated by the elements of the set
  $A=\set{u,v,w,x,y}$, where $u=c$, $v=ac$, $w=ca$, $x=ab$ and
  $y=baba$. It is shown in \cite{CRRT} that there exists no automatic
  structure for $S$ with respect to the generating set $A$. Hence, in
  light of the previous example, a semigroup may be hyperbolic but not
  automatic with respect to a certain choice of generators. It is also
  shown in \cite{CRRT} that finitely generated subsemigroups of free
  semigroups are automatic, so there is some other choice of
  generators for $S$ which supports an automatic structure.
\end{example}

\begin{question}
  Is every hyperbolic semigroup automatic?
\end{question}

We have yet to give an example of a semigroup which is not
hyperbolic. We do so in the next section.

\section{Geometry}

As we mentioned in the introduction, the connection between algebraic
structure and geometry of Cayley graphs is more tenuous for
semigroups than for groups. Nevertheless we will show that a weak form
of the thin triangle condition holds for hyperbolic semigroups. The
triangles in question will not necessarily be geodesic triangles.

Recall that the Cayley graph $\G=\G(S)$ corresponding to a choice of
generators $\Si^+\to S$ for a semigroup $S$ has vertices $S\cup
\set{*}$ where $*=1$ if $S$ has a unit, and otherwise $*$ is a symbol
not in $S$.  For each $a\in \Si$ there is a directed edge labeled $a$
from $*$ to $\ovr a$. Likewise there is a directed edge labeled $a$
from each $s\in S$ to $s\ovr a$. 

The length of the shortest undirected path joining two vertices
defines a metric on the vertices of $\G$. Extend this metric to all of
$\G$ by making edges isometric to the unit interval. For any set
$X\subset \G$, $B_r(X)$ denotes the ball of radius $r$ around
$X$. 

Although it makes sense to speak of paths between any two points in $\G$, 
% ajd 
as in the definition of the metric above, from now on %ajd
we consider only directed paths between vertices. 
For each vertex $p$ and
word $w\in \Sigma^+$ there is a unique path $\gamma_w$ from $p$ to
$p\ovr w$ with label $w$.  We
ignore the distinction between words and paths when the initial point
of the path is understood or irrelevant. Thus we may speak of a path $w$.

Let $D=(u,v,w)$ be a triple of words such that $\ovr{uv}=\ovr w$. For
any vertex $p$ of $\G$, the paths $\gamma_u$, $\gamma_v$ and
$\gamma_w$ beginning at $p$, $p\ovr u$ and $p$ respectively form the
sides of a triangle in $\G$. We shall call $D$ a triangle with
sides $u$, $v$, and $w$. When we talk about the distance between
points on the sides of the triangle, we refer to the distance in $\G$
between points on $\gamma_u$, $\gamma_v$ and $\gamma_w$ for any choice
of $p$. In contrast to the situation for Cayley diagrams of groups,
the choice of $p$ can affect the shape of the triangle.

For any language $L\subset \Sigma^+$, a minimal word of $L$ is
one which has minimum length among all words of $L$ representing the
same element of $S$. A triangle $D=(u,v,w)$ is an $L$--triangle if
$u,v,w\in L$. If in addition $u,v,w$ are minimal words of $L$, then
$D$ is a minimal $L$--triangle. If $R$ is a combing for $S$, then
there is a one to one correspondence between $R$--triangles $(u,v,w)$
and elements $u\#v\#w^r$ in the multiplication table determined by $R$. 

\begin{theorem}\label{thin} If the semigroup $S$ is hyperbolic with respect
  to a choice of generators $\Sigma^+\to S$ and combing $R$, then
  there exists a constant $\delta$ such that for all minimal $R$--triangles
  $D=(u,v,w)$ each point on the side $w$ is a distance at
  most $\delta$ from some point on the union of the other two sides.
\end{theorem}

\begin{proof}
  Let $T$ be the multiplication table corresponding to $R$, and pick a
  context--free grammar $\mathcal G$ in Chomsky normal form for $T$.
  The productions in $\mathcal G$ are of the form
  $X\to YZ$ or $X\to a$ where $X,Y,Z\in N$, the set of nonterminals,
  and 
 %$a\in \Sigma$. 
  $a\in \Sigma_\#$. %ajd 
  Without loss of generality discard nonterminals
  which do not occur in a derivation of any word in $T$. Thus
  for some constant $k$ every nonterminal derives a word $w$
  of length $\abs w \le k$. We show that $\delta=2k$ suffices.

  Fix a minimal $R$--triangle $(u,v,w)$ and
  derivation $S\tto u\#v\#w^r$. Let $q$ be any point
  on the side $w$. Some letter $a$ in $w$ labels an edge of $\Gamma$ 
  containing $q$ or having $q$ as an adjacent vertex. Since $\mathcal
  G$ is in Chomsky normal form, 
  the derivation must include a nonterminal $Z$ which derives that $a$. More
  precisely $w=w_1aw_3$, and the derivation has the form $S\tto \alpha
  Z \beta \to \alpha a \beta \tto u\#v\#w_3^raw_1^r$ where $\alpha\tto
  u\#v\#w_3^r$, $Z\to a$, and $\beta\tto w_1^r$. 

  Adjust the choice of $Z$ so that it derives a subword $w_2^r$ of
  $w^r$ of maximal length containing $a$. Replacing $w_2^r$ by a
word of length no more than $k$ derivable from $Z$ yields an element
$u\#v\#t^r\in T$ and a corresponding triangle $(u,v,t)$. As $\ovr t =
  \ovr {uv} = \ovr w$, minimality of our original 
triangle implies $\abs w \le \abs t$.  Hence $\abs{w_2}\le k$.

$Z$ must occur in the
  derivation as a result of the application of a production
  $X\to YZ$ or $X\to ZY$. The latter case would contradict the
  maximality of $w_2$ as $X$ would derive a subword of $w^r$. Thus
  $X\to YZ$, and $X$ derives a subword 
  containing one or two $\#$'s. Consider the first possibility,
$$
S\tto \alpha X \beta \to \alpha YZ \beta \tto u\#v_1v_2\#w_3^rw_2^rw_1^r
$$
with $Y\tto v_2\#w_3^r$ and $Z\tto w_2^r$. In particular $X\tto
v_2\#w_3^rw_2^r$, and substituting a word $x$ of length at most $k$
derivable from $X$ yields $u\#v_1xw_1^r\in T$. Consequently $x$ has
the form $y\#z^r$, and there is a triangle 
$(u,v_1y, w_1z)$ as in the lefthand side of
Figure~\ref{thintriangle}. We conclude that 
the distance from $q$ to the side $v$ is at most 
$\abs{w_2} + \abs x \le 2k$.
  \begin{figure}[hbt] 
  \centering 
  \psfrag{u}[cc]{$u$}
  \psfrag{u1}[cc]{$u_1$} 
  \psfrag{u2}[cc]{$u_2$}
  \psfrag{v}[lc]{$v$}
  \psfrag{t}[cc]{$t$}
  \psfrag{v1}[lc]{$v_1$}
  \psfrag{v2}[lc]{$v_2$}
  \psfrag{w1}[cc]{$w_1$} 
  \psfrag{w2}[cc]{$w_2$}
  \psfrag{w3}[cc]{$w_3$} 
  \psfrag{p}[cc]{$p$}
  \psfrag{pu}[cc]{$p\ovr u$}
  \psfrag{pw}[cc]{$p\ovr w$}
  \psfrag{z}[cc]{$z$} 
  \psfrag{y}[cc]{$y$}
  \psfrag{q}[cb]{$q$}
  \includegraphics{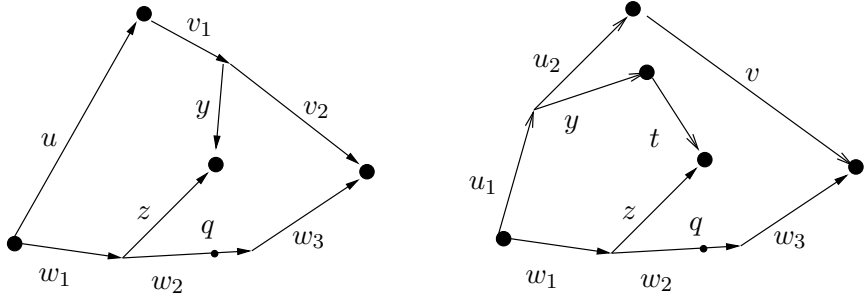}
  \caption{Bounding the distance.\label{thintriangle}}
  \end{figure} 
The second possibility is handled by a similar argument.  The
  subword of $u\#v\#w^r$ derived from $X$ contains two
  $\#$'s, and 
$$
S\tto \alpha X \beta \to \alpha YZ \beta \tto u_1u_2\#v\#w_3^rw_2^rw_1^r
$$
with $Y\tto u_2\#v\#w_3^r$ and $Z\tto w_2^r$. The nonterminal $X$
derives a word 
%$y\#t\#z$ 
 $y\#t\#z^r$ %ajd
of length at most $k$. See the righthand side
of Figure~\ref{thintriangle}. 
\end{proof}

In a cancellative semigroup the argument above works for points on all
sides of minimal $R$--triangles. Thus we have the following theorem,
whose proof is left to the reader.

\begin{theorem}\label{cancellative} If the cancellative
  semigroup $S$ is hyperbolic with respect to a choice of generators
  $A^+\to S$ and combing $R$, then 
  for all minimal $R$--triangles the distance from a point on one side
  of the triangle to the union of the other two sides is uniformly bounded.
\end{theorem}

\begin{corollary}\label{group} A group $G$ is hyperbolic according to
  Definition~\ref{hyperbolic} if and only if it is hyperbolic as a group.  
\end{corollary}
\begin{proof} 
  Suppose the choice of generators $A^+\to G$ and combing $R$ make $G$
  hyperbolic according to Definition~\ref{hyperbolic}. Let $R'$ be the
  the subset of minimal words of $R$. By Theorem~\ref{cancellative}
  the distance from a point on one side of any $R'$--triangle to the
  union of the other two sides is uniformly bounded. As $R'$ is a
  combing (not necessarily regular) of $G$, \cite[Theorem~4]{RG2}
  insures that $G$ is hyperbolic as a group.

  For the converse suppose $G$ is hyperbolic as a
  group. By~\cite[Theorem~1]{RG2} there is a choice of generators
  $\Sigma^*\to G$ and combing $R$ such that the corresponding
  multiplication table is context--free. By Theorem~\ref{monoid} $G$
  is hyperbolic as a semigroup.
\end{proof}

To demonstrate hyperbolicity for a semigroup we need only find a
combing for which the multiplication table is context--free. It seems
more difficult to show that a semigroup is not hyperbolic because we
must survey all possible combings. It helps to have a geometric
condition which is satisfied by all hyperbolic semigroups and which
does not depend on the combing. We start with a lemma modeled
on~\cite[Lemma~1.5 of Chapter~3]{Co} and
\cite[Lemma~9]{RG2}. 

Write $\log(x)$ for $\log_2(x)$, and recall that a
minimal word of a combing $R$ for a semigroup $S$ is one which is of
minimum length among all words in $R$ representing the same element of $S$.

\begin{lemma}\label{nhd} Fix a choice of generators $\Sigma^+\to S$
  and combing $R$ for the hyperbolic semigroup $S$. There are
  constants $k_1, k_2$ such that if $v_0$ and $w$ are paths in $\G$
  from $p$ to $q$ and $v_0$ is a minimal word of $R$ with
  $\ovr{v_0}=\ovr w$, then every point on $v_0$ is a distance at most
  $k_1 + k_2\log(\abs w)$ from $w$.
\end{lemma}

\begin{proof}
  Let $p_0$ be a point on $v_0$. If $\abs w \le 2$, then the distance from
  $p_0$ to $w$ is at most $\abs{v_0} \le k_1$ where $k_1$ is  the
  maximum length of the finitely 
  many minimal words in $R$ representing the images in $S$ of
  words of length 1 or 2 in  $\Sigma^+$. Otherwise estimate the
  distance by constructing a sequence of minimal $R$--triangles as in
  Figure~\ref{triangles}.
  \begin{figure}[ht] 
  \psfrag{p}[cc]{$p$} 
  \psfrag{h}[cc]{$q$} 
  \psfrag{w}[cc]{$w$}
  \psfrag{v0}[cc]{$v_0$} 
  \psfrag{v1}[cc]{$v_1$} 
  \psfrag{v2}[cc]{$v_2$}
  \psfrag{p0}[cc]{$p_0$}
  \psfrag{p1}[cc]{$p_1$} 
  \psfrag{p2}[bc]{$p_2$}
  \includegraphics{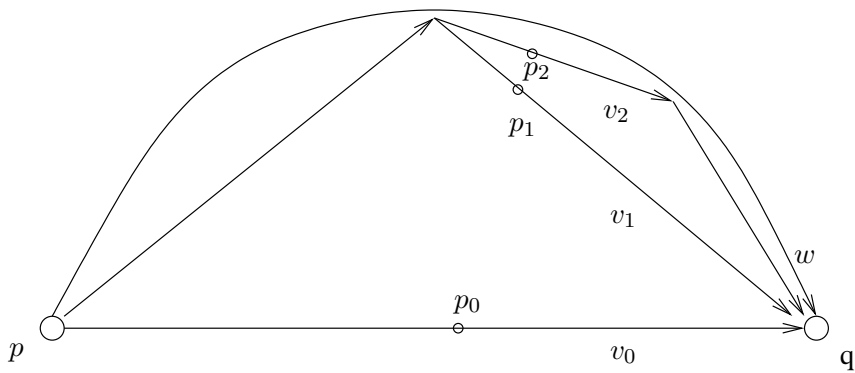}
  \caption{A sequence of minimal $R$--triangles.\label{triangles}}
  \end{figure}

  Start with an $R$--triangle $D_0$ whose base is $v_0$ and whose third
  vertex is a vertex $r$ of $\G(S)$ on the path $w$ and as close to the
  middle of $w$ as possible. Let $w=w_1w_2$ be the corresponding
  decomposition of $w$, and take minimal $R$--words for $w_1$ and
  $w_2$ as the other two sides of $D_0$.

  By Theorem~\ref{thin} there is a side
  $v_1$ of $D_0$ distinct from $v_0$ and containing a point $p_1$ with
  $d(p_1,p_0)\le\delta$. If $v_1$ subtends a segment of $w$ of
  length greater than $2$, construct triangle $D_1$ with base $v_1$ in
  the same way $D_0$ was constructed. Continue until reaching a
  triangle $D_m$ with point $p_{m+1}$ on a side subtending a segment
  of $w$ of length at most $2$.

  To see that the construction of triangles terminates consider the
  sequence of numbers defined by $b_0=\abs w$ and $b_{k+1}=(1 +
  b_k)/2$. It is straightforward to show that $b_k\le 1+\abs w/2^k$.
  By construction the base of triangle $D_k$ subtends a segment of $w$
  of length at most $b_k$. If $b_k\le 4$, then the other two sides of
  $D_k$ subtend segments of length at most $2$ and the construction
  halts.  The estimate above shows that $k\ge \log(\abs w) -1 $ forces
  $b_k\le 4$. As we started with triangle $D_0$, there are no more than
  $k+1=\log(\abs w)$ triangles, and the distance from $p_0$ to $w$ is
  at most $k_1+\delta\log(\abs w)$.
\end{proof}

Here is the geometric condition mentioned above.

\begin{theorem}\label{intersect} Let $A^+\to S$ be a choice of
  generators for the hyperbolic semigroup $S$. There is a constant $k$
with the following property. For any $n\ge 2$ if there is a path of
length at most $n$ from vertex $*$ to vertex $p$ in $\G(S)$, then the
intersection of the balls of radius $k\log n$ around all paths of
length at most $n$ from $*$ to $p$ contains a path from $*$ to $p$.
\end{theorem}
\begin{proof} 
Let $v$ be a minimal combing word for $p$. By Lemma~\ref{nhd} $v$ is in the
ball of radius at most $k_1 + k_2\log n$ around each path of length at
most $n$ from $*$ to $p$. As $n\ge 2$, $k=k_1+k_2$ suffices.
\end{proof}

\begin{example} The semigroup $S=\set{a, b \mid ab=ba}$ is not
  hyperbolic. Choose generators $A=\set{a,b}$ corresponding to the given
  presentation. Consider the paths $w=a^nb^n$ and $v=b^na^n$ from $*$
  to $p=\ovr{a^nb^n}$. By inspection of $\G(S)$ any path of length $n$
  starting at $*$ ends a distance at least $n/2$ from $w$ or from
  $v$.  
  %Since paths from $*$ to $q$ have length $2n$, it follows that
  %the smallest integer $d$ for which $B_d(w)\cap B_d(v)$ 
  %contains a path from $*$ to $\ovr{a^nb^n}$ is at least $n/2$. But if
  %$S$ were hyperbolic, $d$ would be at most $k\log (2n)$. %
  It follows that
  the smallest integer $d$ for which $B_d(w)\cap B_d(v)$ 
  contains a path from $*$ to $p$ is at least $n/2$. But,
  since paths from $*$ to $p$ have length $2n$, if
  $S$ is hyperbolic then  $d$ is  at most $k\log (2n)$.
  %ajd

\end{example}

Corresponding to each choice of generators for a hyperbolic group
there is a canonical combing, namely the set of geodesic
paths. Geodesic paths in $\G(S)$ are not in general
directed and so do not yield combings of $\G(S)$. Nevertheless we can
show that minimal words for the same element of $S$ from different
combings are not too far apart. Notice that if $S$ is hyperbolic with
respect to the choices of generators $A_i^+\to S$ and combings $R_i$
for $i=1,2$, then it is hyperbolic with respect to $(A_1\cup A_2)^+\to
S$ and either combing. Thus when comparing combings we may arrange
things so that both combings are with respect to the same choice of
generators. The following theorem is an immediate consequence of
Lemma~\ref{nhd}. 

\begin{theorem}\label{ft}
  Suppose $S$ is hyperbolic with respect to the choice of generators
  $A\to S$ and combings $R_1$, $R_2$. There are a constants $k_1,k_2$ such
  that if $w_1\in R_1$ and $w_2\in R_2$ are minimal combing paths of
  length at most $n$ for the same element, then each is in the ball of
  radius $k_1 + k_2\log n$ around the other.
\end{theorem} 

Applying Theorem~\ref{ft} to the case $R_1=R_2$, we obtain a weak
analog of the fellow traveler property satisfied by geodesics in hyperbolic
groups. 
  
\section{Word Problems of Semigroups}

The standard word problem of a group $G$ with respect to a choice of
generators $\Sigma^*\to G$ is the language of all words representing
the identity in $G$.  The following definition makes sense for semigroups.

\begin{definition}\label{sg} The word problem for a semigroup $S$ with
  respect to a choice of generators $\Sigma^+\to S$ is the language of
  all words $w\#v^r\in \Sigma_{\#}^*$ such that $w,v\in \Sigma^+$ and
  $\ovr w = \ovr v$ in $S$.  \end{definition}

Recall that a family of languages is called a cone if it contains at
least one non--empty language and is closed under rational transductions.
 
\begin{theorem}\label{cone} Let $\mathcal C$ be a cone.
\begin{enumerate} 
\item \label{p1} If the word problem of $S$ is in $\mathcal C$ for one
  choice of generators, then it is for every choice.
\item\label{p2} If the word problem of $S$ is in $\mathcal C$ for one
  choice of generators, then so is the word problem of every finitely
generated subsemigroup of $S$. 
\end{enumerate} 
\end{theorem} 
\begin{proof} We prove both parts at the same time. 
  Let $f:\Sigma^+\to S$ be a choice of generators for $S$ and
  $g:\Delta^+\to S_0\subset S$ a choice of generators for the
  subsemigroup $S_0$. Denote by $W$ and $W_0$ the corresponding word
  problems, and suppose $W\in \mathcal C$. We must show $W_0\in
  \mathcal C$.

  Pick a
  homomorphism $h:\Delta^+\to\Sigma^+$ such that $f\circ h=g$.
  By Lemma~\ref{reverse-lemma}
  $\rho=(h)(\#,\#)(\tau_h)$ is a rational transduction. It is apparent
  that $\rho$ is a partial function with domain $\Delta^+\#\Delta^+$
  and that $\rho(x\#y^r)=h(x)\#h(y)^r$. Now $x\#y^r\in W_0$ if and
  only if $f(h(x))=g(x)=g(y)=f(h(y))=f(\tau_h(y^r)^r)$. Hence
  $x\#y^r\in W_0$ if and only if $h(x)\#\tau_h(y^r)=\rho(x\#y^r)\in
  W$. In other words, $W_0=\rho\inv(W)$ whence $W_0\in \mathcal C$.
\end{proof}

\begin{theorem}\label{grouplanguage} Let $\mathcal C$ be a cone and
  $G$ a group. The standard word problem for $G$ is in $\mathcal C$
  if and only if the semigroup word problem for $G$ is in
  $\mathcal C$.
\end{theorem}
\begin{proof}

By the preceding theorem we are free to make a convenient choice of
generators. 
Choose a surjective homomorphism $\Sigma^*\to G$ such that the
restriction to $\Sigma^+$ is onto and $\Sigma$ has formal inverses. 
The latter condition means that 
$\Sigma=\set{ a, a\inv, b, b\inv, \ldots }$ with $\ovr {a
\inv}=\ovr a\inv$ etc. There is an induced involutory map $w\to w\inv$
of $\Sigma^*$ to itself such that $(a_1\cdots a_n)\inv=a_n\inv\cdots
a_1\inv$ and $\ovr{w\inv}=\ovr w\inv$. 

Let $V=\set{w\mid w\in\Sigma^*, \ovr w = 1}$ and $W=\set{u\#v^r\mid
u,v\in \Sigma^+, \ovr u = \ovr v }$ be the standard and semigroup word
problems for $G$ respectively.

Choose $w_1\in\Sigma^+$ with $\ovr{w_1}=1$, and observe that $V=\e +
\set{w \mid w\#w_1^r\in W}$. It follows that $V$ is the image of $W$
under the union of the rational transduction which maps $\Sigma_\#^*$ to
$\e$ and the rational transduction $\rho=(I)(\#,\e)(w_1^r,\e)$. Here
$I:\Sigma^*\to \Sigma^*$ is the identity map. As rational
transductions are closed under union, $V$ must be in $\mathcal C$ if
$W$ is. 

For the converse suppose $V\in\mathcal C$, and let $h:\Sigma_\#^*\to
\Sigma^*$ be the homomorphism 
which is the identity on $\Sigma$ and maps $\#$ to $\e$. As inverse
homomorphisms are rational transductions, $h\inv(V)$ is in $\mathcal
C$ and consequently so is $V_1=h\inv(V)\cap (\Sigma^+\#\Sigma^+) =
\set{x\#y\mid x,y\in\Sigma^+, \ovr{xy}=1}$. 

Define $\tau:\Sigma^+\to\Sigma^+$ by $\tau(w)=(w\inv)^r$;
$\tau=\bigl(\sum_{a\in\Sigma}(a,a\inv)\bigr)^+$ is a rational
transduction. Consequently $\rho=(I)(\#,\#)(\tau)$ is a rational
transduction from $\Sigma_\#^*$ to itself. 

We claim that $W$ is the image of $V_1$ under $\rho$. If $x\#y\in
V_1$, then $\rho(x\#y)=x\#(y\inv)^r$; and $\ovr{xy}=1$ implies
$x\#(y\inv)^r\in W$. Conversely if $u\#v^r\in W$, then $\ovr u =\ovr v$
implies $u\#v\inv\in V_1$; and $\rho(u\#v\inv)=u\#v^r$.
\end{proof}

It is not hard to show that The word problem of $S$ is regular if and
  only if $S$ is finite.  Muller and Schupp's
  characterization of groups with context--free word
  problem~\cite{MS1} (see also~\cite{Du}) suggests
  the following question.

\begin{question} What is the structure of semigroups with context--free
 word problem?
\end{question}

\end{document}